\documentclass[preprint,11pt]{article}
\usepackage{amsmath,amssymb}
\usepackage{multirow}
\usepackage{enumerate}
\usepackage{array}
\usepackage{algorithm}
\usepackage{algorithmic}

\newtheorem{theo}{Theorem}[section]

\newtheorem{lemm}[theo]{Lemma}
\newtheorem{prop}[theo]{Proposition}

\usepackage{xcolor}

\begin{document}

\begin{center}
{\Large \bf Isomorphism and isotopism classes of filiform Lie algebras of dimension up to seven}
\end{center}

\begin{center}
{\large \em O. J. Falc\'on$^1$,\ R. M. Falc\'on$^2$\, J. N\'u\~nez$^{3}$}
\vspace{0.5cm}

{\small $^{1, 3}$ Faculty of Mathematics, University of Seville, Spain.\\
$^2$School of Building Engineering, University of Seville, Spain.\\
E-mail: {\em $^1$oscfalgan@yahoo.es,\ $^2$rafalgan@us.es,\ $^3$jnvaldes@us.es}}
\end{center}

\vspace{0.5cm}

\noindent {\large \bf Abstract.} Since the introduction of the concept of isotopism of algebras by Albert in 1942, a prolific literature on the subject has been developed for distinct types of algebras. Nevertheless, there barely exists any result on the problem of distributing Lie algebras into isotopism classes. The current paper is a first step to deal with such a problem. Specifically, we define a new series of isotopism invariants and we determine explicitly the distribution into isotopism classes of $n$-dimensional filiform Lie algebras, for $n\leq 7$. We also deal with the distribution of such algebras into isomorphism classes, for which we confirm some known results and we prove that there exist $p+8$ isomorphism classes of seven-dimensional filiform Lie algebras over the finite field $\mathbb{F}_p$ if $p\neq 2$.

\vspace{0.5cm}


\section{Introduction.}

At present, the problem of the classification of Lie algebras in general is
still unsolved. According to Levi's theorem, there exist three types of Lie algebras: simple, solvable and semi-direct sum of semi-simple and solvable Lie algebras. Semi-simple Lie algebras include simple Lie algebras and solvable Lie algebras contain nilpotent Lie algebras. In turn, filiform Lie algebras constitute the most structured subset of nilpotent Lie algebras and have a large number of applications in Applied Mathematics, Engineering and Physics \cite{Georgi1999, Gilmore2005}. They were introduced formally by Vergne \cite{Vergne1966, Vergne1970} in the late 1960s, although Umlauf had already used them as an example in his thesis \cite{Umlauf1891}. The distribution into isomorphism classes of $n$-dimensional filiform Lie algebras over the complex field is known for $n\leq 12$ \cite{Boza2003, Fedriani1997}, even if it is only known for nilpotent Lie algebras over the complex field of dimension $n\leq 7$ \cite{Ancochea1989}. More recently, some authors have dealt with the classification of $n$-dimensional nilpotent Lie algebras over finite fields $\mathbb{F}_p=\mathbb{Z}/p\mathbb{Z}$. Specifically, Schneider \cite{Schneider2005} obtained the number of isomorphism classes over the finite field $\mathbb{F}_2$, for $n\leq 9$, and over $\mathbb{F}_3$ and $\mathbb{F}_5$, for $n\leq 7$. The classification of six-dimensional nilpotent Lie algebras over a field of characteristic distinct of two was determined by de Graaf \cite{Graaf2007} and over any arbitrary field by Cical\`o, de Graaf and Schneider \cite{Cicalo2012}. With respect to the classification of filiform Lie algebras over $\mathbb{F}_p$, Schneider obtained that

\begin{theo}\label{theoSch}\cite{Schneider2005} There exist six six-dimensional filiform Lie algebras over $\mathbb{F}_2$ and five over $\mathbb{F}_3$ and $\mathbb{F}_5$; $15$ seven-dimensional filiform Lie algebras over $\mathbb{F}_2$, $11$ over $\mathbb{F}_3$ and $13$ over $\mathbb{F}_5$; $47$ eight-dimensional filiform Lie algebras over $\mathbb{F}_2$ and $124$ nine-dimensional filiform Lie algebras over $\mathbb{F}_2$. \hfill $\Box$
\end{theo}

\vspace{0.5cm}

The main goal of this paper is to step forward in the distribution of filiform Lie algebras over any field, not only into isomorphism classes, which is the usual criterion, but also into isotopism classes, which is, at present, the first contribution on this subject. Isotopisms of filiform Lie algebras constitute a generalization of isomorphisms that can be used to gather together non-isomorphic algebras. The concept of isotopism of algebras was introduced by Albert \cite{Albert1942} and developed later by Bruck \cite{Bruck1944a}, who introduced the concept of {\em isotopically simple algebra} as a simple algebra such that all their isotopic algebras are simple. He delved further into the study of isotopisms of division algebras and simple algebras. Isotopisms have been used since then to study and classify distinct algebraic structures like Jordan algebras \cite{Petersson1969}, alternative algebras \cite{Babikov1997} or division algebras \cite{Schwarz2010}. Nevertheless, there barely exists any result about isotopisms of Lie algebras, apart from the following two results of Albert and Bruck:

\begin{lemm} \cite{Albert1942} \label{lemm_iso0} A principal isotope $\mathfrak{g}$ of a Lie algebra $\mathfrak{h}$ with respect to an isotopism $(f,g,\mathrm{Id})$ is a Lie algebra if and only if the following two conditions are verified.
\begin{enumerate}[i.]
\item $[f(u),g(v)]=-[f(v),g(u)]$, for all $u,v\in\mathfrak{g}$.
\item $[f([f(u),g(v)]),g(w)] - [f([f(u),g(w)]),g(v)] - [f(u),g([f(v),g(w)])]= 0$, for all $u,v,w\in\mathfrak{g}$. \hfill $\Box$
\end{enumerate}
\end{lemm}

\newpage

\begin{theo} \cite{Bruck1944a} \label{theo_iso0} It is verified that:
\begin{enumerate}[i.]
\item The Lie algebra of order $n(n-1)/2$, consisting of all skew-symmetric matrices, over any subfield of the field of all reals, under the multiplication $A\circ B=AB-BA$, is isotopically simple.
\item The  Lie  algebra of order $n(n-1)$, consisting of all skew-hermitian matrices in any field $R(i)$ (where $R$ is a subfield of the reals and $i^2= —1$), under  the  multiplication $A\circ B=AB-BA$, is an isotopically simple algebra over $R$. \hfill $\Box$
\end{enumerate}
\end{theo}

\vspace{0.2cm}

More recently, Jim\'enez-Gestal and P\'erez-Iquierdo \cite{Jimenez2008} have studied the relationship that exists between the isotopisms of a finite-dimensional real division algebra and the Lie algebra of its ternary derivations. On the other hand, Allison et al. \cite{Allison2009, Allison2012} have recently studied isotopes of a class of graded Lie algebras called Lie tori, but the notion of isotopism that they use is quite different from the conventional one.

\vspace{0.2cm}

The structure of the paper is as follows: Section 2 is devoted to some preliminary concepts and results on Lie algebras that we use throughout the paper. We define a new series of isotopism invariants, which constitute the initial point of the distribution of six- and seven-dimensional filiform Lie algebras into isomorphism and isotopism classes that we determine, respectively, in Sections 3 and 4. In them, we also analyze those nonlinear equations that relate the structure constants of two isomorphic filiform Lie algebras of such dimensions with the entries of a matrix of an isomorphism between them. Such an analysis yields Theorem \ref{theoSch} and also the existence of $p+8$ isomorphism classes of seven-dimensional filiform Lie algebras, over $\mathbb{F}_p$, for any prime $p>2$.

\section{Preliminaries.}

We start with some basic concepts and results on Lie algebras that we use throughout the paper. For a more comprehensive review on the subject, the reader can consult \cite{Varadarajan1984}. An $n$-dimensional {\em Lie algebra} $\mathfrak{g}$ over a field $\mathbb{K}$ is an $n$-dimensional vector space over $\mathbb{K}$ endowed with a second inner law, named {\em bracket product}, that is bilinear and anti-commutative and that satisfies the {\em Jacobi identity}
\begin{equation}\label{ij}
J(u,v,w)=[u,[v,w]]+[v,[w,u]]+[w,[u,v]]=0, \mbox{ for all }
u,v,w\in \mathfrak{g}.
\end{equation}
Given a basis $\{e_1,\ldots,e_n\}$ of the Lie algebra $\mathfrak{g}$, the {\em structure constants} of $\mathfrak{g}$ are the numbers $c_{ij}^k\in \mathbb{K}$ such that
\begin{equation}\label{filbrack2}
[e_i,e_j] = \sum_{k=1}^n c_{ij}^k e_k, \hspace{0.5cm} \mbox{for} \hspace{0.5cm} 1 \leq i < j \leq n.
\end{equation}
The {\em centralizer} of a subset $\mathfrak{h}$ of $\mathfrak{g}$ is the set $\mathrm{Cen}_{\mathfrak{g}}(\mathfrak{h})=\{u\in\mathfrak{g}\mid\, [u,v]=0, \text{ for all } v\in \mathfrak{h}\}$. An {\em ideal} of $\mathfrak{g}$ is any vector subspace $\mathfrak{h}\subseteq \mathfrak{g}$ such that $[\mathfrak{h},\mathfrak{g}]\subseteq \mathfrak{h}$. It is called {\em abelian} if $[\mathfrak{h},\mathfrak{h}]=0$. The {\em center} of $\mathfrak{g}$ is the ideal $Z(\mathfrak{g})=\mathrm{Cen}_{\mathfrak{g}}(\mathfrak{g})$. Given an $m$-dimensional ideal $\mathfrak{h}$ of $\mathfrak{g}$,  the {\em quotient Lie algebra} of $\mathfrak{g}$ over $\mathfrak{h}$ is the quotient space $\mathfrak{g}/\mathfrak{h}$. It has structure of $n-m$ dimensional Lie algebra with the product
\begin{equation}\label{quot}
[u+\mathfrak{h},v+\mathfrak{h}]=[u,v]+\mathfrak{h}, \text{ for all } u,v\in \mathfrak{g}.
\end{equation}

Two Lie algebras $\mathfrak{g}$ and $\mathfrak{g}'$ are {\em isotopic} if there exist three regular linear transformations $f,g$ and $h$ from $\mathfrak{g}$ to $\mathfrak{g}'$ such that
\begin{equation}\label{isot}
[f(u),g(v)]=h([u,v]), \text{ for all } u,v\in\mathfrak{g}.
\end{equation}
It is denoted as $\mathfrak{g}\simeq\mathfrak{g}'$ and it is said that $\mathfrak{g}'$ is an {\em isotope} of $\mathfrak{g}$. The tuple $(f,g,h)$ is called an {\em isotopism} of Lie algebras. If $h$ is the identity transformation $\mathrm{Id}$, then the isotopism is said to be {\em principal}. If $f=g=h$, then it is an isomorphism of Lie algebras. It is denoted as $\mathfrak{g}\cong\mathfrak{g}'$ and the isomorphism $(f,f,f)$ is denoted only as $f$. Given two isomorphic Lie algebras, $\mathfrak{g}$ and $\mathfrak{g}'$, an isomorphism $f$ between them and an ideal $\mathfrak{h}$ of $\mathfrak{g}$, the vector subspace $f(\mathfrak{h})$ is an ideal of $\mathfrak{g}'$ and the quotient Lie algebras $\mathfrak{g}/\mathfrak{h}$ and $\mathfrak{g}'/f(\mathfrak{h})$ are also isomorphic. If $\mathfrak{h}=Z(\mathfrak{g})$, then $f(\mathfrak{h})=Z(\mathfrak{g}')$. Hereafter, $\mathfrak{g}^{(i)}$ represents the quotient Lie algebra $\mathfrak{g}^{(i-1)}/Z(\mathfrak{g}^{(i-1)})$, for all natural $i\in\mathbb{N}$, where $\mathfrak{g}^{(0)}=\mathfrak{g}$.

To be isotopic and isomorphic are equivalence relations on Lie algebras, which can then be partitioned into {\em isotopism} and {\em isomorphism classes}. Since isotopisms of Lie algebras constitute a generalization of isomorphism of such algebras, the problem of distributing Lie algebras into isotopism classes is reduced to that of gathering together isomorphism classes with isotopic equivalence class representatives. Even if there exist in the literature distinct isomorphism invariants that are traditionally used in the classification of Lie algebras into isomorphism classes, there does not exist any previous study about isotopism invariants. In the current paper, we define a sequence of such invariants, which also are, evidently, isomorphism invariants for Lie algebras. Specifically, given an $n$-dimensional Lie algebra $\mathfrak{g}$ and a natural $m\leq n$, we define
\begin{equation}\label{isotinv}
d_m(\mathfrak{g})=\max\{\dim \mathrm{Cen}_{\mathfrak{g}}(\mathfrak{h})\colon\, \mathfrak{h} \text{ is an } m\text{-dimensional ideal of } \mathfrak{g}\}.
\end{equation}
The next result proves that the sequence $d(\mathfrak{g})=\{d_1(\mathfrak{g}),\ldots,d_n(\mathfrak{g})\}$ is an isotopism invariant.

\begin{lemm} \label{lmm_inv} Given two isotopic $n$-dimensional Lie algebras, $\mathfrak{g}$ and $\mathfrak{g}'$, and a natural $m\leq n$, it is $d_m(\mathfrak{g})=d_m(\mathfrak{g}')$.
\end{lemm}

{\bf Proof} Let $(f,g,h)$ be an isotopism between $\mathfrak{g}$ and $\mathfrak{g}'$. Given $m\leq n$, let $\mathfrak{h}$ be an $m$-dimensional ideal of $\mathfrak{g}$ such that $d_m(\mathfrak{g})=\dim \mathrm{Cen}_{\mathfrak{g}}(\mathfrak{h})$. Given $x\in\mathrm{Cen}_{\mathfrak{g}}(\mathfrak{h})$, it is $[f(x),g(\mathfrak{h})]=h([x,\mathfrak{h}])=h(0)=0$ and hence, $f(\mathrm{Cen}_{\mathfrak{g}}(\mathfrak{h}))\subseteq \mathrm{Cen}_{\mathfrak{g}'}(g(\mathfrak{h}))$. Since $f$ and $g$ are regular transformations, it is $\dim g(\mathfrak{h}) = m$ and $d_m(\mathfrak{g})\leq d_m(\mathfrak{g}')$. The equality follows from a similar reasoning that makes use of the isotopism $(f^{-1},g^{-1},h^{-1})$ between $\mathfrak{g}'$ and $\mathfrak{g}$. \hfill $\Box$

\vspace{0.5cm}

The {\em lower central series} of a Lie algebra $\mathfrak{g}$ is defined as the series of ideals
\begin{equation}
{\small\mathcal{C}^1(\mathfrak{g})=
\mathfrak{g} \ \supseteq
\mathcal{C}^2(\mathfrak{g})=[\mathfrak{g},\mathfrak{g}] \ \supseteq \dots\ \supseteq \mathcal{C}^k(\mathfrak{g})=[\mathcal{C}^{k-1}(\mathfrak{g}),
\mathfrak{g}] \ \supseteq \dots}
\end{equation}
If there exists a natural $m$ such that $\mathcal{C}^m(\mathfrak{g}) \equiv 0$, then $\mathfrak{g}$ is called {\em nilpotent}. Nilpotency is preserved by isomorphisms of Lie algebras, because $f(\mathcal{C}^i(\mathfrak{g}))=\mathcal{C}^i(\mathfrak{g}')$, for every isomorphism $f$ between $\mathfrak{g}$ and an isomorphic Lie algebra $\mathfrak{g}'$. The {\em nil-index} of $\mathfrak{g}$ is the smallest natural $p$ such that $\mathcal{C}^{p+1}(\mathfrak{g})\equiv 0$. Since $\mathcal{C}^p(\mathfrak{g})\subseteq Z(\mathfrak{g})$, the center of a non-trivial nilpotent Lie algebra is non-trivial. The {\em type} of $\mathfrak{g}$ is defined as the sequence
\begin{equation}\label{type}
\{\dim \mathfrak{g}/\mathcal{C}^2(\mathfrak{g}),\, \dim \mathcal{C}^2(\mathfrak{g})/\mathcal{C}^3(\mathfrak{g}),\ldots,\dim \mathcal{C}^{p-1}(\mathfrak{g})/\mathcal{C}^p(\mathfrak{g})\}.
\end{equation}
According to Engel's theorem, the Lie algebra $\mathfrak{g}$ is isomorphic to a nilpotent Lie algebra of basis $\{e_1,\ldots,e_n\}$, such that
\begin{equation}\label{engel}
\begin{array}{c}
\begin{cases}
[e_1,e_2]=0,\\
[e_1,e_i]=\epsilon_{i-1} e_{i-1}, \text{ for all } i\in\{3,\ldots,n\},
\end{cases}
\end{array}
\end{equation}
where $\epsilon_2,\ldots,\epsilon_{n-1}\in\{0,1\}$. The Lie algebra $\mathfrak{g}$ is said to be {\em filiform} if $\dim \, \mathcal{C}^k(\mathfrak{g}) = n-k$, for all $k\in \{2,\ldots, n\}$. Filiformity is also preserved by isomorphisms of Lie algebras. Further, the next two numbers were defined by Echarte et al. \cite{Echarte1996, Echarte1996a} as isomorphism invariants for filiform Lie algebras.
\begin{equation}\label{z1}
z_1(\mathfrak{g})=\max \{k\in\mathbb{N}\mid\, \mathrm{Cen}_{\mathfrak{g}}(\mathcal{C}^{n-k+2}(\mathfrak{g}))\supset \mathcal{C}^2(\mathfrak{g})\}.
\end{equation}
\begin{equation}\label{z2}
z_2(\mathfrak{g})=\max \{k\in\mathbb{N}\mid\, \mathcal{C}^{n-k+1}(\mathfrak{g})\text{ is abelian}\}.
\end{equation}
The only bidimensional filiform Lie algebra is the abelian. For $n\geq 3$, the type of an $n$-dimensional filiform Lie algebra $\mathfrak{g}$ is $\{2,1,\ldots,1\}$. Since the center of $\mathfrak{g}$ is one-dimensional, it is verified that $d_1(\mathfrak{g})=n$ and $d_n(\mathfrak{g})=1$. A basis $\{e_1,\ldots,e_n\}$ of $\mathfrak{g}$ is said to be {\em compatible with respect to its lower central series} if
\begin{equation}\label{lcs}
\mathfrak{g}^2 =\langle e_2,\ldots ,e_{n-1}\rangle,\, \mathfrak{g}^3 =\langle e_2,\ldots ,e_{n-2} \rangle,\ldots, \mathfrak{g}^{n-1}=\langle e_2\rangle,
\,\, \mathfrak{g}^n=0.
\end{equation}
Hence,
\begin{equation}\label{lcs1}
Z(\mathfrak{g})=\langle e_2\rangle
\end{equation}
and
\begin{equation}\label{lcs2}
[e_i,e_j] \in \langle e_2,\ldots ,e_{i-1}\rangle, \text{ whenever } 3\leq i<j\leq n.
\end{equation}
Since $[e_3,u]\in\langle e_2\rangle$ and $J(e_1,e_3,u)=0$, for all  $u\in\mathfrak{g}$, the fact of being $e_2\in Z(\mathfrak{g})$ implies that
\begin{equation}\label{lcs3}
[e_3,u]=0, \text{ for all } u\in \mathfrak{g}^2.
\end{equation}
From (\ref{lcs}), if $n\geq 4$, then $[e_4,e_n]=c_{4n}^2e_2+c_{4n}^3e_3$. Thus, $0=J(e_1,e_4,e_n)=[e_1,c_{4n}^2e_2+c_{4n}^3e_3]-[e_4,e_{n-1}]$ and hence,
\begin{equation}\label{lcs4}
[e_4,e_{n-1}]\in \langle e_2\rangle.
\end{equation}

If the base field of an $n$-dimensional filiform Lie algebra $\mathfrak{g}$ has characteristic zero, then Vergne \cite{Vergne1966, Vergne1970} proved the existence of an {\em adapted basis} that facilitates the distribution of filiform Lie algebras into isomorphism classes. It can be considered, without loss of generality, to be compatible with respect to its lower central series, verifying Conditions (\ref{lcs1}-\ref{lcs4}), and can be chosen in such a way that the following conditions are also satisfied:
\begin{equation}\label{lcs5}
\begin{array}{c}
\begin{cases}
[e_1,e_{i+1}]=e_i, \text{ for all } i\in\{3,\ldots,n-1\},\\
[e_3,e_n]=0.
\end{cases}
\end{array}
\end{equation}
If the only brackets that are distinct of zero in an adapted basis are those of the form $[e_1,e_{i+1}]=e_i$, then the corresponding filiform Lie algebra is called {\em model}. The $n$-dimensional model algebra is the only filiform Lie algebra having the isomorphism invariant $z_1$ equal to $n$ and hence, it is not isomorphic to any other algebra of the same dimension.

If the base field has characteristic distinct of zero, then it is not always possible to find an adapted basis of a filiform Lie algebra $\mathfrak{g}$. Due to it and similarly to the definition of a {\em nilpotent basis} proposed by Schneider in \cite{Schneider2005}, we define a {\em filiform basis} of $\mathfrak{g}$ as a compatible basis $\{e_1,\ldots,e_n\}$ with respect to its lower central series such that, given $i\in\{3,\ldots,n\}$, it is $[e_1,e_i]=e_{i-1}$ or $[e_i,e_n]=e_{i-1}$. The following results hold.

\begin{lemm}\label{lemm_fb} Regardless of the base field, every finite-dimensional filiform Lie algebra has a filiform basis.
\end{lemm}

{\bf Proof} Let $\mathfrak{g}$ be a filiform Lie algebra over a field $\mathbb{K}$. Since $\mathfrak{g}$ is nilpotent, there exists a basis $\beta=\{e_1,\ldots,e_n\}$ satisfying Condition (\ref{engel}) and being compatible with respect to the lower central series of $\mathfrak{g}$. Let $i\in\{3,\ldots,n\}$ be such that $[e_1,e_i]=0$. From (\ref{lcs}) and (\ref{lcs2}), there exists $a_i\in \mathbb{K}\setminus\{0\}$ such that $[a_ie_i,e_n]=e_{i-1}$. We can suppose all these numbers $a_i$ to be equal to $1$. Otherwise, we replace each $e_i$ by $a_ie_i$ in the basis $\beta$. \hfill $\Box$

\vspace{0.5cm}

\begin{prop}\label{prop_iso} The $n$-dimensional model algebra is not isotopic to any other filiform Lie algebra of the same dimension.
\end{prop}

{\bf Proof} Given a non-model $n$-dimensional filiform Lie algebra $\mathfrak{g}$ of filiform basis $\{e_1,\ldots,e_n\}$, there must exists a non-zero bracket $[e_i,e_j]$, such that $i\neq 1\neq j$. It implies that $d_{n-1}(\mathfrak{g})<n-1$. However, if $\mathfrak{h}$ is the $n$-dimensional model algebra, then $d_m(\mathfrak{h})=n-1$, for all $m\in\{2,\ldots,n-1\}$. The result follows then from Lemma \ref{lmm_inv}. \hfill $\Box$

\vspace{0.5cm}

The only isomorphism class of the set of filiform Lie algebras of dimension $n\leq 4$ corresponds to the model algebra of such a dimension. It coincides, therefore, with the only isotopism class of such a set. For $n=5$, there exist two isomorphism classes of filiform Lie algebras: the model algebra and that having an adapted basis satisfying the bracket $[e_4,e_5]=e_2$. Both classes determine indeed distinct isotopism classes, because the isotopism invariant $d_4$ is equal to $4$ for the model algebra, but it is equal to $2$ for the algebras of the second class. For higher dimensions, the distribution into isomorphism and isotopism classes requires a more detailed study of the corresponding constant structures. We analyze them in the next sections, where we make use of the following lemma for the corresponding distribution into isomorphism classes.

\begin{lemm} \label{lemm0} Given a filiform Lie algebra $\mathfrak{g}$ of filiform basis $\{e_1,\ldots,e_n\}$ and an isomorphism $f$ between $\mathfrak{g}$ and a filiform Lie algebra $\mathfrak{g'}$, it is verified that $f(e_i)\in \mathcal{C}^{n-i+1}(\mathfrak{g}')$, for all $i\in\{2,\ldots,n-1\}$.
\end{lemm}

{\bf Proof} The result follows straightforward from (\ref{lcs}) and (\ref{lcs2}). \hfill $\Box$

\vspace{0.5cm}

Given a filiform basis $\{e'_1,\ldots,e'_n\}$ of $\mathfrak{g'}$, let $F=(f_{ij})$ be the regular square matrix that defines $f$, that is to say, such that $f(e_j)=\sum_{i=1}^n f_{ij}e'_i$, for all natural $j\leq n$. Since $e_2\in Z(\mathfrak{g})$ and $\mathfrak{g}^2 =\langle e_2,\ldots ,e_{n-1}\rangle$, the numbers $f_{21}$ and $f_{2n}$ do not have any
influence on the isomorphism. We can suppose that $f_{21}=f_{2n}=0$ and hence, Lemma \ref{lemm0} yields that
\begin{equation}\label{F}
{\scriptsize F = \left(\begin{array} {cccccccc}
f_ {11} & 0 & 0 & \ldots & 0 & f_{1n}  \\
0 & f_{22} & f_{23} & \ldots & f_{2(n-1)} & 0 \\
f_{31} & 0 & f_{33} & \ldots & f_{3(n-1)} & f_{3n} \\
\vdots & \vdots & \vdots & \ddots & \vdots & \vdots \\ f_{(n-1)1} & 0 & 0
& \ldots & f_{(n-1)(n-1)} & f_{(n-1)n} \\ f_{n1} & 0 &  0 & \ldots & 0 & f_{nn} \end{array} \right)}
\end{equation}

\section{Classification of six-dimensional filiform Lie algebras.}

The distribution into isomorphism classes of six-dimensional filiform Lie algebras is known over any field \cite{Cicalo2012, Graaf2007, Schneider2005}. Nevertheless, it is interesting to deal again with such a distribution, independently of the known results, in order to take advantage of such a small dimension to expose in detail the algebraic methodology related to the new isotopism invariants. Since the reasoning is analogous for higher dimensions, we avoid in the following section to show explicitly all the steps, because of their extensive and tedious computation. A case study of possible six-dimensional filiform bases allows us to assure the next result.

\begin{lemm}\label{f6} Given a six-dimensional filiform Lie algebra $\mathfrak{g}$ over a field $\mathbb{K}$, there exist three numbers $a,b,c\in \mathbb{K}$ and an adapted basis of $\mathfrak{g}$ such that
\begin{equation}
\mathfrak{g}\cong\mathfrak{g}^6_{abc}\equiv\begin{cases}
[e_1,e_{i+1}]=e_i, \text{ for all } i>1,\\
[e_4,e_5]=a e_2,\\
[e_4,e_6]=b e_2 + a e_3,\\
[e_5,e_6]=c e_2 + b e_3 + a e_4.
\end{cases}
\end{equation}
\hfill$\Box$
\end{lemm}

\vspace{0.2cm}

The isotopism invariants of the previous algebras determine a preliminary classification for their distribution into isotopism and isomorphism classes.
\begin{equation}\label{D6}
d(\mathfrak{g}^6_{abc})=\begin{cases}
\{6,5,5,5,5,1\}, \text{ if } a=b=c=0,\\
\{6,5,3,2,2,1\}, \text{ if } a=b=0\neq c,\\
\{6,5,4,4,2,1\}, \text{ if } a=0\neq b,\\
\{6,5,3,2,2,1\}, \text{ if } a\neq 0.
\end{cases}
\end{equation}
Lemma \ref{lmm_inv} yields then the existence of at least four isotopism (isomorphism) classes of six-dimensional filiform Lie algebras over any field. Since the model algebra $\mathfrak{g}^6_{000}$ constitutes an isotopism (isomorphism) class by itself, we focus on the rest of cases. We start with the distribution into isomorphism classes.

\begin{prop}\label{propF6_2} Given $c,C\in \mathbb{K}\backslash \{0\}$, it is $\mathfrak{g}^6_{00c}\cong\mathfrak{g}^6_{00C}$.
\end{prop}

{\bf Proof} It is enough to consider the isomorphism $f$ between them, defined as $\mathbb{K}$-linear extension of
$$f(e_i)=\begin{cases}
e_1, \text{ if } i=1,\\
\frac cC e_i, \text{ otherwise}.
\end{cases}$$ \hfill $\Box$

\vspace{0.5cm}

The remaining cases require an analysis of the entries of the matrix $F=(f_{ij})$ of an isomorphism between two six-dimensional isomorphic filiform Lie algebras over a field $\mathbb{K}$, $\mathfrak{g}^6_{abc}$ and $\mathfrak{g}^6_{ABC}$, of respective adapted bases $\{e_1,\ldots,e_6\}$ and $\{e'_1,\ldots,e'_6\}$. Its form has been exposed in (\ref{F}). Further, since $-f_{33}f_{16}e'_2=[f(e_3),f(e_6)]=f([e_3,e_6])=f(0)=0$ and $F$ is regular, it must be $f_{16}=0$. The central columns of $F$ can be determined by using the entries of the first and last columns.

\begin{lemm} \label{lemmF6} It is verified that
\begin{enumerate}[a)]
\item $f_{55}=f_{11}f_{66}$.
\item $f_{44}=(f_{11} - A f_{61}) f_{55}$.
\item $f_{33}=(f_{11}-A f_{61}) f_{44}$.
\item $f_{22}=f_{11}f_{33}$.
\item $f_{45}=(f_{11}- A f_{61}) f_{56} +A f_{51}f_{66}$.
\item $f_{34}=(f_{11} - A f_{61}) f_{45} - B f_{55}f_{61}$.
\item $f_{23}=f_{11}f_{34} - (Af_{51}+Bf_{61})f_{44}$.
\item $f_{35}=(f_{11}- A f_{61})f_{46}+A f_{41}f_{66}+B(f_{51}f_{66}-f_{56}f_{61})$.
\item $f_{24}=f_{11}f_{35}-(Af_{51}+Bf_{61})f_{45} + (Af_{41}-Cf_{61})f_{55}$.
\item $f_{25}=f_{11}f_{36}+A(f_{41}f_{56}-f_{46}f_{51})+ B(f_{41}f_{66}-f_{46}f_{61}) + C(f_{51}f_{66}-f_{56}f_{61})$.
\end{enumerate}
\end{lemm}

{\bf Proof} The result follows from the fact that
$f(e_i)=f([e_1,e_{i+1}])=[f(e_1),$ $f(e_{i+1})]$, for all $i\in\{2,\ldots,n-1\}$. \hfill $\Box$

\vspace{0.5cm}

Since $F$ is a regular matrix, we have then that
\begin{equation}\label{detF6}
f_{11}f_{66}(f_{11}-A f_{61})\neq 0.
\end{equation}
This last condition together with the rest of brackets related to Expression (\ref{isot}) yield the rest of constraints for the entries of $F$. Specifically, it must be
\begin{equation}\label{eq1}
a f_{11} = A (f_{66}+a f_{61}),
\end{equation}
\begin{equation}\label{eq2}
b (f_{11}-A f_{61})^2= B (f_{66}+a f_{61})=0,
\end{equation}
\begin{equation}\label{F6_cf}
\begin{array}{c}f_{11}f_{66}(cf_{11}(f_{11}-A f_{61})^2
+bB(Af_{61}-2f_{11})f_{61}+ 2aAf_{41}-aCf_{61}\\+2Af_{46}-Cf_{66})-Af_{56}^2(f_{11}-A f_{61})
-A^2(af_{51} + 2f_{56})f_{51}f_{66}=0.
\end{array}
\end{equation}
The following result holds.

\begin{prop}\label{propF6_3} Given two filiform Lie algebras, $\mathfrak{g}^6_{0bc}$ and $\mathfrak{g}^6_{0BC}$, over a field $\mathbb{K}$, such that $b\neq 0\neq B$, both algebras are isomorphic if $\mathbb{K}$ has characteristic distinct of two. Otherwise, the algebras $\mathfrak{g}^6_{010}$ and $\mathfrak{g}^6_{011}$ are not isomorphic.
\end{prop}

{\bf Proof} After imposing $a=A=0$ in Equations  (\ref{detF6} - \ref{F6_cf}), it results that
$$\begin{cases}
bf_{11}^2=Bf_{66},\\
cf_{11}^4-2B^2f_{61}f_{66}=Cf_{11}f_{66}.
\end{cases}$$
Any isomorphism between $\mathfrak{g}^6_{0bc}$ and $\mathfrak{g}^6_{0BC}$ verifies then that $f_{66}=bf_{11}^2/B$. If the characteristic of the base field is two, then we have that $f_{11}=f_{66}=1$ and $c=C$. Hence, $\mathfrak{g}^6_{010}$ is not isomorphic to $\mathfrak{g}^6_{011}$. Now, if the characteristic is distinct or two, then it is enough to consider $f_{61}=(Bcf_{11}^2-bCf_{11})/2B^2b$ to obtain the required isomorphism. \hfill $\Box$

\vspace{0.5cm}

Equation (\ref{F6_cf}) can be used to determine an entry of the matrix $F$ that does not appear in Equations (\ref{detF6} - \ref{eq2}). In particular, if $a\neq 0\neq A$, then we isolate the variable $f_{46}$ in Equation (\ref{F6_cf}) if the characteristic of the base field $\mathbb{K}$ is distinct of two, or the variable $f_{56}$, otherwise. The following result is then deduced from Condition (\ref{detF6}) and the isotopism invariants (\ref{D6}).

\newpage

\begin{theo}\label{theoF6_1} It is verified that
\begin{enumerate}[a)]
\item Any two filiform Lie algebras, $\mathfrak{g}^6_{a0c}$ and $\mathfrak{g}^6_{A0C}$, over a field $\mathbb{K}$, such that $a\neq 0\neq A$, are isomorphic.
\item Any two filiform Lie algebras, $\mathfrak{g}^6_{abc}$ and $\mathfrak{g}^6_{ABC}$, over a field $\mathbb{K}$, such that $a\neq 0\neq A$ and $b\neq 0\neq B$, are isomorphic.
\item None of the algebras of (a) is isomorphic to one of (b). \hfill $\Box$
\end{enumerate}
\end{theo}

\vspace{0.2cm}

The previous results establish the following distribution into isomorphism classes of the six-dimensional filiform Lie algebras.

\begin{theo}\label{theoF6_2} There exist six isomorphism classes of six-dimensional filiform Lie algebras over a field of characteristic two:
$$\mathfrak{g}^6_{000}, \mathfrak{g}^6_{001}, \mathfrak{g}^6_{010}, \mathfrak{g}^6_{011}, \mathfrak{g}^6_{100} \text{ and } \mathfrak{g}^6_{110}.$$ If the base field has characteristic distinct of two, then there exist five isomorphism classes, which correspond to the previous ones, keeping in
mind that now, $\mathfrak{g}^6_{010}\cong\mathfrak{g}^6_{011}$. \hfill $\Box$
\end{theo}

\vspace{0.2cm}

We finish the current section with the distribution of six-dimensional filiform Lie algebras into isotopism classes.

\begin{prop}\label{propi6} There exist five isotopism classes of six-dimensional filiform Lie algebras over any field:
$$\mathfrak{g}^6_{000},\mathfrak{g}^6_{001}, \mathfrak{g}^6_{010}, \mathfrak{g}^6_{100} \text{ and }\mathfrak{g}^6_{110}.$$
\end{prop}

{\bf Proof} From Theorem \ref{theoF6_2} and the isotopism invariants (\ref{D6}), there exist at least four isotopism classes of six-dimensional filiform Lie algebras. We observe that $\mathfrak{g}^6_{011}\simeq \mathfrak{g}^6_{010}$ by the isotopism $(f,f,h)$, where $f$ and $h$ are defined as the respective linear extension of
$$f(e_4)=e_4-e_3 \text{ and } f(e_i)=e_i, \text{ if } i\neq 4,$$
$$h(e_3)=e_3-e_2 \text{ and } h(e_i)=e_i, \text{ if } i\neq 3.$$
On the other hand, a simple analysis on possible isotopisms between the Lie algebras $\mathfrak{g}^6_{100}$ and $\mathfrak{g}^6_{110}$ determines that these two algebras are not isotopic, whatever the base field is. Therefore, there exist exactly five isotopism classes of six-dimensional filiform Lie algebras over any field. \hfill $\Box$

\section{Classification of seven-dimensional filiform Lie algebras.}

A case study of possible seven-dimensional filiform bases gives us the next result.

\begin{lemm}\label{f7} Given a seven-dimensional filiform Lie algebra $\mathfrak{g}$ over a field $\mathbb{K}$ of characteristic distinct of two, there exist four numbers $a,b,c,d\in \mathbb{K}$ and an adapted basis of $\mathfrak{g}$ such that
\begin{equation}\label{fb7a}
\mathfrak{g}\cong\mathfrak{g}^7_{abcd}\equiv\begin{cases}
[e_1,e_{i+1}]=e_i, \text{ for all } i>1,\\
[e_4,e_7]=a e_2,\\
[e_5,e_6]=b e_2,\\
[e_5,e_7]=c e_2 + (a+b) e_3,\\
[e_6,e_7]=d e_2 + c e_3 + (a+b) e_4.
\end{cases}
\end{equation}
If the base field $\mathbb{K}$ has characteristic two, then there exists a filiform basis of $\mathfrak{g}$ of one of the following three non-isomorphic types:
\begin{enumerate}
\item[a)] {\bf Type 1:} An adapted basis satisfying (\ref{fb7a}).
\item[b)] {\bf Type 2:} A filiform basis such that
\begin{equation}\label{fb7b}
\mathfrak{g}\cong\mathfrak{g}^7_a\equiv\begin{cases}
[e_1,e_3]=[e_4,e_6]=[e_5,e_7]=e_2,\\
[e_4,e_7]=[e_5,e_6]=e_3,\\
[e_1,e_5]=e_4,\\
[e_1,e_6]=e_5,\\
[e_1,e_7]=e_6,\\
[e_6,e_7]=e_3 + a e_4, \text{ for some } a\in \{0,1\}.\\
\end{cases}
\end{equation}
\item[c)] {\bf Type 3:} A filiform basis such that
\begin{equation}\label{fb7c}
\mathfrak{g}\cong\mathfrak{h}^7_a\equiv\begin{cases}
[e_3,e_7]=[e_4,e_6]=e_2,\\
[e_1,e_4]=[e_5,e_6]=e_3,\\
[e_5,e_7]=e_4,\\
[e_6,e_7]=e_5,\\
[e_1,e_7]=e_6,\\
[e_4,e_7]=ae_2, \text{ for some } a\in \{0,1\}.\\
\end{cases}
\end{equation}\hfill$\Box$
\end{enumerate}
\end{lemm}

\vspace{0.2cm}

The isotopism invariants related to the previous algebras are
\begin{equation}\label{D7}
d(\mathfrak{g}^7_{abcd})=\begin{cases}
\{7,6,6,6,6,6,1\}, \text{ if } a=b=c=d=0,\\
\{7,6,6,6,5,4,1\}, \text{ if } a=b=c=0\neq d,\\
\{7,6,6,5,5,3,1\}, \text{ if } a=b=0\neq c,\\
\{7,6,6,4,3,3,1\}, \text{ if } a=0\neq b,\\
\{7,6,5,5,5,2,1\}, \text{ if } a\neq 0=b,\\
\{7,6,5,4,3,2,1\}, \text{ if } a\neq 0\neq b.\\
\end{cases}
\end{equation}
\begin{equation}\label{D7a}
d(\mathfrak{g}^7_0)=d(\mathfrak{g}^7_1)=\{7,6,4,4,2,2,1\}.
\end{equation}
\begin{equation}\label{D7b}
d(\mathfrak{h}^7_0)=d(\mathfrak{h}^7_1)=\{7,6,5,4,3,2,1\}.
\end{equation}
Hence, the algebras $\mathfrak{g}^7_0$ and $\mathfrak{g}^7_1$ are not isotopic to any other type algebra of Lemma \ref{f7}. Indeed, they determine two distinct isotopism classes, because their quotient Lie algebras ${\mathfrak{g}^7_0}^{(2)}$ and ${\mathfrak{g}^7_1}^{(2)}$ have distinct isotopism invariants. Specifically, $d_3({\mathfrak{g}^7_0}^{(2)})=4\neq 3 = d_3({\mathfrak{g}^7_1}^{(2)})$. Further, the algebras $\mathfrak{h}^7_0$ and $\mathfrak{h}^7_1$ are not isotopic to any algebra of first type, because the isotopism invariants of their quotient Lie algebras are distinct. Specifically,
\begin{equation}\label{D7q}
d({\mathfrak{g}^7_{abcd}}^{(1)})=\begin{cases}
\{6,5,5,5,5,1\}, \text{ if } a+b=c=0,\\
\{6,5,5,4,3,1\}, \text{ if } a+b=0\neq c,\\
\{6,5,4,4,2,1\}, \text{ if } a+b\neq 0.\\
\end{cases}
\end{equation}
\begin{equation}\label{D7b2}
d({\mathfrak{h}^7_0}^{(1)})=d({\mathfrak{h}^7_1}^{(1)})=\{6,5,4,3,2,1\}.
\end{equation}
Nevertheless, it is $\mathfrak{h}^7_0\simeq\mathfrak{h}^7_1$, by the principal isotopism $(f,f,\mathrm{Id})$, where $f$ is defined by $\mathbb{K}$-linear extension from $f(e_4)=e_3+e_4$ and $f(e_i)=e_i$, for all $i\neq 4$.

We focus then our study on seven-dimensional filiform Lie algebras with adapted bases. We start with their distribution into isomorphism classes. Let $\mathfrak{g}^7_{abcd}$ and $\mathfrak{g}^7_{ABCD}$ be two seven-dimensional isomorphic filiform Lie algebras over a field $\mathbb{K}$ and let $F=(f_{ij})$ be the matrix related to an isomorphism $f$ between them. Analogously to six-dimensional filiform Lie algebras, the  central columns of $F$ are determined by its first and last columns. In particular,
\begin{equation}\label{f7diag}
f_{ii}=f_{11}^{7-i}f_{77}, \text{ for all } i\in\{2,\ldots,6\}.
\end{equation}
Since $F$ is a regular matrix, it is
\begin{equation}\label{detF7}
f_{11}f_{77}\neq 0
\end{equation}
and the following equations are then deduced from the definition of isomorphism
\begin{equation}\label{eq3}
af_{11}^2=Af_{77},
\end{equation}
\begin{equation}\label{eq3a}
bf_{11}^2=Bf_{77},
\end{equation}
\begin{equation} \label{F7_cf1}
cf_{11}^4 =2(A+B)^2f_{71}f_{77}+Cf_{11}f_{77},
\end{equation}
\begin{equation} \label{F7_cf2}
\begin{array}{c}
df_{11}^5f_{77}-(3A+2B)cf_{11}^3f_{71}f_{77}+aA^2f_{11}f_{71}^2f_{77}+\\
(bA^2+aAB+bAB)f_{11}f_{71}^2f_{77}-2(AC+BC)f_{72}f_{77}^2-\\
Bf_{11}f_{67}^2+2Bf_{11}f_{57}f_{77}=Df_{11}f_{77}^2.
\end{array}
\end{equation}
The following results hold.

\begin{prop}\label{propF7_2}
Two non-model algebras $\mathfrak{g}^7_{00cd}$ and $\mathfrak{g}^7_{00CD}$ are isomorphic if and only if one of the following conditions is verified
\begin{enumerate}[a)]
\item $c=C=0$ and $d\neq 0\neq D$.
\item $d=D=0$ and $c\neq 0\neq C$.
\item The numbers $c,C,d$ and $D$ are all of them distinct of zero.
\end{enumerate}
\end{prop}

{\bf Proof} Once we impose $a=A=b=B=0$ in Equations (\ref{detF7} - \ref{F7_cf2}), it results that $cf_{11}^3=Cf_{77}$ and $df_{11}^4=Df_{77}$. Hence, it is enough to impose $f_{77}=df_{11}^4/D$ if we are in the first case or $f_{77}=cf_{11}^3$ if we are in the second one. In the third
case, it must be $f_{77}=df_{11}^4/D=cf_{11}^3/C$ and thus, we impose $f_{11}=cD/Cd$. \hfill $\Box$

\vspace{0.5cm}

\begin{prop}\label{propF7_3} Let $\mathfrak{g}^7_{0bcd}$ and $\mathfrak{g}^7_{0BCD}$ be two non-model
algebras such that $b\neq 0\neq B$. They are isomorphic if and
only if one of the following assertions is verified
\begin{enumerate}
\item[a)] The characteristic of the base field is distinct of two.
\item[b)] The characteristic of the base field is two and $c=C$.
\end{enumerate}
\end{prop}

{\bf Proof} Once we impose $a=A=0$ in Equations (\ref{detF7} - \ref{F7_cf2}), we obtain that
$$\begin{cases}
bf_{11}^2=Bf_{77},\\
cf_{11}^4=2B^2f_{71}f_{77} + Cf_{11}f_{77},\\
df_{11}^5f_{77}=2B(cf_{11}^3f_{71}f_{77}+Cf_{71}f_{77}^2-f_{11}f_{57}f_{77}) +Bf_{11}f_{67}^2+Df_{11}f_{77}^2.
\end{cases}$$
Hence, any isomorphism between $\mathfrak{g}^7_{0bcd}$ and
$\mathfrak{g}^7_{0BCD}$ must verify that $f_{77}=bf_{11}^2/B$. If
the base field has characteristic distinct of two, then it is enough to isolate $f_{71}$ and $f_{57}$ in the second and third equations, respectively. Otherwise, it is $f_{11}=f_{77}=b=B=1$ and therefore, $c=C$ and $d=f_{67}+D$. Thus, it is enough to take $f_{67}=d-D$ to obtain our isomorphism. \hfill $\Box$

\vspace{0.5cm}

\begin{prop}\label{propF7_4} It is verified that
\begin{enumerate}
\item[a)] If the base field $\mathbb{K}$ has characteristic two, then $\mathfrak{g}^7_{1000}\cong\mathfrak{g}^7_{1001}$, but none pair of filiform Lie algebras of the set $\{\mathfrak{g}^7_{1000},\mathfrak{g}^7_{1010}, \mathfrak{g}^7_{1011}\}$ are isomorphic.
\item[b)] If $\mathbb{K}$ has characteristic distinct of two, then every filiform Lie algebra $\mathfrak{g}^7_{a0cd}$ such that $a\neq 0$ and $4ad=5c^2$ is isomorphic to $\mathfrak{g}^7_{1000}$. If $4ad\neq 5c^2$ and $4ad-5c^2$ is a perfect square in $\mathbb{K}$, then it is isomorphic to $\mathfrak{g}^7_{1001}$.
\item[c)] If $\mathbb{K}$ is an algebraically closed field of characteristic distinct of two, then every filiform Lie algebra $\mathfrak{g}^7_{a0cd}$ such that $a\neq 0$ and $4ad\neq 5c^2$ is isomorphic to $\mathfrak{g}^7_{1001}$.
\item[d)] Given a prime $p\neq 2$ and a non-perfect square $q$ in $\mathbb{K}=\mathbb{F}_p$, every filiform Lie algebra $\mathfrak{g}^7_{a0cd}$ such that $a\neq 0$, $4ad\neq 5c^2$ and $4ad-5c^2$ is a non-perfect square in $\mathbb{F}_p$, is isomorphic to $\mathfrak{g}^7_{100q}$.
\end{enumerate}
\end{prop}

{\bf Proof} Let $\mathfrak{g}^7_{a0cd}$ and $\mathfrak{g}^7_{A0CD}$ be two filiform Lie algebras such that $a\neq 0\neq A$. Once we impose $b=B=0$ in Equations (\ref{detF7}) - (\ref{F7_cf2}), we obtain that
$$\begin{cases}
af_{11}^2=Af_{77},\\
cf_{11}^4=2A^2f_{71}f_{77} + Cf_{11}f_{77},\\
df_{11}^5=3Acf_{11}^3f_{71}-aA^2f_{11}f_{71}^2+2ACf_{71}f_{77}+Df_{11}f_{77}
\end{cases}$$
Any isomorphism between $\mathfrak{g}^7_{a0cd}$ and
$\mathfrak{g}^7_{A0CD}$ must then verify that $f_{77}=af_{11}^2/A$. If the base field has characteristic two, then $f_{11}=f_{77}=a=A=1$ and the only conditions to study are $c=C$ and $d=cf_{71}-f_{71}+D$. If $c=1$, then $d=D$; otherwise, we take $f_{71}=D-d$ and we have proved the first assertion.

If the characteristic of the base field is distinct of two, then $f_{71}=(cf_{11}^4-Cf_{11}f_{77})/2A^2f_{77}=(Acf_{11}^2-aCf_{11})/2aA^2$
and thus,
$$(4ad-5c^2)A^2f_{11}^2=(4AD-5C^2)a^2.$$
Hence, the algebras $\mathfrak{g}^7_{a0cd}$ and $\mathfrak{g}^7_{A0CD}$ are isomorphic if and only if one of the following conditions is satisfied
\begin{enumerate}
\item[i.] $4ad=5c^2$ and $4AD=5C^2$.
\item[ii.] $4ad\neq 5c^2$, $4AD\neq 5C^2$ and $(4AD-5C^2)/(4ad-5c^2)$ is a perfect square in $\mathbb{K}$.
\end{enumerate}
Assertions (b) and (c) follow then immediately from (i) and (ii). Finally, let $p$ be a prime distinct of two and let $q$ be a non-perfect square in $\mathbb{K}=\mathbb{F}_p$. Every perfect square $r$ in $\mathbb{F}_p$ is uniquely related to a non-perfect square $s$ in $\mathbb{F}_p$ such that $r=q/s$. Such a relation is $1$-$1$, because $\mathbb{F}_p$ is a finite field that contains exactly $(p-1)/2$ perfect squares distinct of zero. Assertion (d) follows then from (ii). \hfill $\Box$

\vspace{0.5cm}

\begin{prop}\label{propF7_5} Let $\mathfrak{g}^7_{abcd}$ and $\mathfrak{g}^7_{ABCD}$ be two non-model algebras such that $a,b,A$ and $B$ are all of them distinct of zero. They are isomorphic if and only $aB=Ab$ and one of the following assertion
is verified
\begin{enumerate}
\item The characteristic of the base field is distinct of two and $a+b\neq 0$.
\item If the characteristic of the base field is two or $a+b=0$, then it must be $c=0=C$ or $c\neq 0\neq C$.
\end{enumerate}
\end{prop}

{\bf Proof} If the characteristic of the base field is distinct of two, then we can isolate $f_{57}$ from Equation (\ref{F7_cf2}). Otherwise, we isolate $f_{67}$. Further, since $A\neq 0$, we have from (\ref{eq3}) that
$f_{77}=af_{11}^2/A=bf_{11}^2/B$ and hence, $aB=Ab$. Equation (\ref{F7_cf1}) yields then that
$$aAcf_{11}^2-a^2Cf_{11}=2(a+b)^2A^2f_{71}.$$
If the characteristic of the base field is distinct of two and $a+b\neq 0$, then a possible isomorphism can be determined once $f_{71}$ is isolated from the previous equation. Otherwise, it is $Acf_{11}=aC$ and hence, there exists an isomorphism between $\mathfrak{g}^7_{abcd}$ and $\mathfrak{g}^7_{A(Ab/a)CD}$ only if $c=C=0$ or $c\neq 0\neq C$. In the last case, it is enough to consider $f_{11}=aC/Ac$. \hfill $\Box$

\vspace{0.5cm}

The previous results establish the distribution into isomorphism classes of seven-dimensional filiform Lie algebras.

\newpage

\begin{theo}\label{theoF7_2} It is verified that
\begin{enumerate}[a)]
\item There exist $15$ isomorphism classes of seven-dimensional filiform Lie algebras over a field of characteristic two:
$$\mathfrak{g}^7_{0000}, \mathfrak{g}^7_{0001}, \mathfrak{g}^7_{0010}, \mathfrak{g}^7_{0011},\mathfrak{g}^7_{0100},\mathfrak{g}^7_{0110}, \mathfrak{g}^7_{1000}, \mathfrak{g}^7_{1010},$$ $$\mathfrak{g}^7_{1011}, \mathfrak{g}^7_{1100}, \mathfrak{g}^7_{1110},\mathfrak{g}^7_0, \mathfrak{g}^7_1, \mathfrak{h}^7_0 \text{ and } \mathfrak{h}^7_1.$$
\item Any seven-dimensional filiform Lie algebra over an algebraically closed field of characteristic distinct of two is isomorphic to one of the following isomorphism classes:
$$\mathfrak{g}^7_{0000}, \mathfrak{g}^7_{0001}, \mathfrak{g}^7_{0010}, \mathfrak{g}^7_{0011},\mathfrak{g}^7_{0100}, \mathfrak{g}^7_{1001}, \mathfrak{g}^7_{1b00}, \text{ and } \mathfrak{g}^7_{1(-1)10},$$
where $b\in\mathbb{K}$.
\item Given a prime $p\neq 2$, there exist $p+8$ isomorphism classes of seven-dimensional filiform Lie algebras over the finite field $\mathbb{F}_p$:
$$\mathfrak{g}^7_{0000}, \mathfrak{g}^7_{0001}, \mathfrak{g}^7_{0010}, \mathfrak{g}^7_{0011},\mathfrak{g}^7_{0100},\mathfrak{g}^7_{1001},\mathfrak{g}^7_{100q},\mathfrak{g}^7_{1b00} \text{ and } \mathfrak{g}^7_{1(-1)10},$$
where $b\in\mathbb{K}$ and $q$ is a non-perfect square of $\mathbb{F}_p$.\hfill $\Box$
\end{enumerate}
\end{theo}

\vspace{0.2cm}

Finally, the distribution of seven-dimensional filiform Lie algebras into isotopism classes is determined in the next result.

\begin{prop}\label{propi7} It is verified that
\begin{enumerate}[a)]
\item There exist 10 isotopism classes of seven-dimensional filiform Lie algebras over a field of characteristic two:
$$\mathfrak{g}^7_{0000}, \mathfrak{g}^7_{0001}, \mathfrak{g}^7_{0010}, \mathfrak{g}^7_{0100},\mathfrak{g}^7_{1000}, \mathfrak{g}^7_{1100}, \mathfrak{g}^7_{1110}, \mathfrak{g}^7_0, \mathfrak{g}^7_1 \text{ and } \mathfrak{h}^7_0.$$
\item There exist eight isotopism classes of seven-dimensional filiform Lie algebras over an algebraically closed field of characteristic distinct of two and also over the finite field $\mathbb{F}_p$, where $p$ is a prime distinct of two:
$$\mathfrak{g}^7_{0000}, \mathfrak{g}^7_{0001}, \mathfrak{g}^7_{0010}, \mathfrak{g}^7_{0100},\mathfrak{g}^7_{1000}, \mathfrak{g}^7_{1100}, \mathfrak{g}^7_{1(-1)00} \text{ and } \mathfrak{g}^7_{1(-1)10}.$$
\end{enumerate}
\end{prop}

{\bf Proof} We have already seen that $\mathfrak{h}^7_0\simeq \mathfrak{h}^7_1$ and that the algebras $\mathfrak{g}^7_0$ and $\mathfrak{g}^7_1$ constitute two distinct isotopism classes by themselves. We determine also the following isotopisms, which are defined as linear extensions from the exposed non-trivial changes of the corresponding bases:
\begin{enumerate}[a)]
\item $\mathfrak{g}^7_{0011}\simeq \mathfrak{g}^7_{0010}$, by the isotopism $(f,f,h)$, where
$$f(e_4)=e_4-e_3\hspace{0.5cm} \text{ and } \hspace{0.5cm} h(e_3)=e_3-e_2.$$
\item $\mathfrak{g}^7_{0110}\simeq \mathfrak{g}^7_{0100}$, by the isotopism $(f,f,h)$, where
$$f(e_1)=e_1+e_7, \hspace{1cm} f(e_7)=e_6+e_7,$$
$$h(e_4)=e_4-e_3, \hspace{1cm} h(e_5)=e_5-e_4 \hspace{0.5cm} \text{ and } \hspace{0.5cm} h(e_6)=e_6+e_5-e_4.$$
\item $\mathfrak{g}^7_{10cd}\simeq \mathfrak{g}^7_{1000}$, for all $c,d\in \mathbb{K}$, by the isotopism $(f,f,h)$, where
$$f(e_5)=e_5 + (c^2-d)e_3, \hspace{1cm} h(e_3)=e_3 - c e_2,$$ $$h(e_4)=e_4 - c e_3 + (c^2-d)e_2 \hspace{0.5cm} \text{ and } \hspace{0.5cm} h(e_5)=e_5 - c e_4.$$
\item $\mathfrak{g}^7_{1b00}\simeq \mathfrak{g}^7_{1100}$, for all $b\in \{2,\ldots,p-2\}$, by the isotopism $(f,f,h)$, where
$$f(e_3)= f_{44}e_3, \hspace{0.5cm} f(e_4)= f_{44}e_4, \hspace{0.5cm} f(e_5)= f_{55}e_5, \hspace{0.5cm} f(e_6)= f_{66}e_6,$$
$$h(e_2)= f_{44}e_2, \hspace{0.5cm} h(e_3)= f_{44}e_3, \hspace{0.5cm} h(e_4)= f_{55}e_4, \hspace{0.5cm} h(e_5)= f_{66}e_5,$$
with
$$f_{44} = \frac{4}{(b+1)^2} f_{66}, \hspace{0.25cm} f_{55} = \frac{2}{b+1} f_{66}, \hspace{0.25cm}
f_{66} = \frac{2b}{b+1}
$$
\end{enumerate}

A simple analysis on possible isotopisms among the Lie algebras $\mathfrak{g}^7_{1100}$, $\mathfrak{g}^7_{1(-1)00}$ and  $\mathfrak{g}^7_{1(-1)10}$ indicates that they determine distinct isotopism classes. The result follows then from Lemma \ref{lmm_inv} and the invariants (\ref{D7})-(\ref{D7b}).  \hfill $\Box$

\section{Final remarks and further work.}

In the current paper we have dealt with the distribution of Lie algebras into isotopism classes. We have defined a series of isotopism invariants of Lie algebras, which we have used to determine the isotopism and isomorphism classes of filiform Lie algebras of dimension $n\leq 7$. We have generalized the study of Schneider \cite{Schneider2005} since we have obtained that there exist $p+8$ isomorphism classes of seven-dimensional filiform Lie algebras over the finite field $\mathbb{F}_p$, for any prime $p\neq 2$. Higher dimensions can be similarly analyzed, although with more extensive case studies. On the other hand, it remains open the problem of distributing into isotopism classes other families of Lie algebras, distinct of the filiform case.

\end{document}